\documentclass[12pt]{amsart}
\usepackage{amsmath,amssymb}
\usepackage{graphicx, verbatim}
\usepackage{enumerate,enumitem}
\usepackage{mathtools}
\usepackage[retainorgcmds]{IEEEtrantools}
\usepackage{youngtab}
\usepackage{xcolor}
\definecolor{mygray}{gray}{.9}
\definecolor{mygray2}{gray}{.8}
\usepackage[pagebackref=true]{hyperref}
\usepackage{cleveref}
\usepackage[margin=1.5in]{geometry}
\usepackage[all]{xy}
\usepackage{tikz}
\usetikzlibrary{shapes,decorations}
\usepackage{aliascnt}
\newtheorem{theorem}{Theorem}[section]

\newaliascnt{lemma}{theorem}
\newtheorem{lemma}[lemma]{Lemma}
\aliascntresetthe{lemma}

\newaliascnt{corollary}{theorem}
\newtheorem{corollary}[corollary]{Corollary}
\aliascntresetthe{corollary}

\newaliascnt{proposition}{theorem}
\newtheorem{proposition}[proposition]{Proposition}
\aliascntresetthe{proposition}

\theoremstyle{definition}
\newtheorem{definition}[theorem]{Definition}

\newtheorem{remark}[theorem]{Remark}
\newtheorem{example}[theorem]{Example}
\newtheorem{question}[theorem]{Question}

\numberwithin{equation}{section}

\DeclareMathOperator{\relint}{{\rm relint}}

\DeclareMathOperator{\Hom}{\rm Hom}

\DeclareMathOperator{\Frac}{Frac}
\newcommand{\sM}{\mathcal{M}}

\newcommand{\blank}{\rule[-1pt]{9pt}{1pt}}
\newcommand{\bm}{\mathfrak{m}}
\DeclareMathOperator{\rank}{rank}
\newcommand{\sC}{\mathcal{C}}
\DeclareMathOperator{\interior}{int}

\begin{document}
\bibliographystyle{alpha}
\newcommand{\br} {{r}}                              
\newcommand{\A} {\mathbb{A}}                           
\newcommand{\PGL} {\Pj\Gl_2(\R)}                           
\newcommand{\RP} {\R\Pj^1}                                 
\newcommand{\C} {{\mathbb C}}                              
\newcommand{\bF}{\mathbb F}
\newcommand{\bR}{\mathbb R}
\newcommand{\bZ}{\mathbb Z}
\newcommand{\cD} {{\mathcal D}}
\newcommand{\cE} {{\mathcal E}}
\newcommand{\cF} {{\mathcal F}}
\newcommand{\cG} {{\mathcal G}}

\newcommand{\cI}{\mathcal I}
\newcommand{\cL} {{\mathcal L}}
\newcommand{\gr} {\mathrm{gr}}
\newcommand{\lev}{\mathrm{lev}}
\newcommand{\inn}{\mathrm{in}}
\newcommand{\R} {{\mathbb R}}                              
\newcommand{\Q}{{\mathbb Q}}                              
\newcommand{\Z} {{\mathbb Z}}                              
\newcommand{\N}{{\mathbb N}}                               
\renewcommand{\k}{{\mathbf k}}                               

\newcommand{\Pj} {{\mathbb P}}                             
\newcommand{\T} {{\mathbb T}}                              
\newcommand{\Sg} {{\mathbb S}}                             
\newcommand{\Gl} {{\rm Gl}}                                
\newcommand{\Sl} {{\rm Sl}}                                
\newcommand{\SL}{{\mathrm SL}}
\newcommand{\fn}{\mathfrak n}

\newcommand{\F}  [2] {\ensuremath{C_{#2}({#1})}}                            
\newcommand{\oM} [1] {\ensuremath{{\mathcal M}_{0,#1}(\R)}}                 
\newcommand{\M} [1] {\ensuremath{{\overline{\mathcal M}}{_{0, #1}(\R)}}}    
\newcommand{\cM} [1] {\ensuremath{{\mathcal M}_{0, #1}}}                    
\newcommand{\CM} [1] {\ensuremath{{\overline{\mathcal M}}_{0, #1}}}         
\newcommand{\bd} [1] {\mathbf{#1}}

\newcommand{\roverm}[1]{\frac{r(#1)}{d}}

\newcommand{\suchthat} {\ \ | \ \ }
\newcommand{\ore} {\ \ {\it or} \ \ }
\newcommand{\oand} {\ \ {\it and} \ \ }

\newcommand{\vw}{\mathbf{w}}
\newcommand{\PPP}[2]{(\Pj^{#1})^{#2}}
\newcommand{\floor}[1]{\lfloor #1 \rfloor}
\newcommand{\ceiling}[1]{\lceil #1 \rceil}
\newcommand{\ceil}[1]{\lceil #1 \rceil}
\def\Proj{\mathrm{Proj} \,}
\def\q{/\!/}
\newcommand{\tableau}[2]{\begin{array}{|c|}\hline #1 \\ \hline #2 \\ \hline \end{array}}
\newcommand{\vd}{\mathbf{d}}
\newcommand{\va}{\mathbf{a}}
\newcommand{\angles}[1]{\langle #1 \rangle}
\newcommand{\oq}{\frac{1}{q}}

\def\O{\mathcal{O}}                               
\def\X{\mathcal{X}}
\def\L{\mathcal{L}}
\def\<{\langle}
\def\>{\rangle}
\def\({\left (}
\def\){\right )}
\def\Vol{\mathrm{Vol}}

\title { 
$F$-splittings of seminormal monoid algebras}
\date\today
\subjclass[2000]{Primary 14P25, Secondary 90C48, 52B11}

\keywords{}

\author{Milena Hering}
\address{Milena Hering, University of Edinburghi and Maxwell Institute for Mathematical Sciences, Edinburgh, UK}
\email{m.hering@ed.ac.uk}

\author{Kevin Tucker}
\address{Kevin Tucker, University of Illinois at Chicago, Chicago, IL, USA }
\email{kftucker@uic.edu}

\begin{abstract}
	We compute a number of invariants of singularities defined via the Frobenius morphism for seminormal affine toric varieties over fields of characteristic $p > 0$. Our main technical tool is a combinatorial description of the potential splittings of iterates of Frobenius for seminormal monoid algebras. This allows us to give an easy formula for the F-splitting ratio of such rings  as well as to compute the ideals stable under the Cartier algebra, including the test ideal.
\end{abstract}

\thanks{We are grateful to the Max Planck Institut Bonn and Mathematical Sciences Research Institute for their hospitality. 
 The first author
   is grateful for support from 
   NSF grant
   DMS \#1001859,  EPSRC First grant EP/K041002/1, a LMS Emmy Noether Fellowship, and EPSRC Fellowship EP/T018836/1. The second author is thankful for support under NSF Grants DMS \#1602070 and \#1707661, and for a
fellowship from the Sloan Foundation.}
   
\maketitle

\section{Introduction}
Let $k$ be a perfect field of characteristic $p$, let $R$ be a finitely generated $k$-algebra, and let $F\colon R\to R$ be the Frobenius morphism. The study of $F$-splittings or $R$-module sections of Frobenius have featured prominently in many diverse fields, with numerous applications in commutative algebra, algebraic geometry, number theory, and representation theory. Our goal in this paper is to describe the set $\Hom_R(R^{1/p^e},R)$ of so-called $p^{-e}$-linear maps and potential $F^e$-splittings when $R$ is a seminormal monoid ring. This then gives us the tool to compute uniformly $F$-compatible ideals, including the test ideal,  and the $F$-splitting dimension and $F$-splitting ratio governing the asymptotic growth rate of the number of iterated Frobenius splittings. 

Seminormal monoids have been studied for example by  Li \cite{Lithesis}, Bruns, Li, R\"omer \cite{BrunsLiRoemer}, and De Stefani, Monta\~no, and N{\'u}\~nez-Betancourt \cite{StefaniMontanoNunezBetancourt} and the goal of this article is to extend their study. While seminormal monoids are more complex than normal monoids, they are still sufficiently tractable to provide a theory that makes it easy to check conjectures, give examples and thus build a testing ground for further exploration. 

Let $S$ be an affine monoid embedded in a lattice $M$,  and let $C$ be the rational cone generated by $S$ in $M\otimes \bR$. We assume $C$ is pointed. When $S$ is seminormal, the points in the interior of a face $D$ are given by the intersection of that face with a sublattice of $M$, denoted $M_D$ (see Definition \ref{def:seminormal}). We say that a face $D$ is relatively unsaturated, or an RUF, if $M_D$ is strictly smaller than the intersection of the lattices of the faces strictly containing $D$.
One of the main points of this article is that 
the RUFs are the key to understanding the properties of the  seminomormal monoid algebra. 
We say that a face is a $p$-face if $p$ divides the index of $M_D$ in $M\cap \<D\>$ and call a RUF that is also a $p$-face a $p$-RUF, see Definition \ref{def:pface}.  Note that a monoid ring is $F$-split if and only if it is seminormal and has no $p$-faces (see Section \ref{sec:basics}  and Theorem \ref{thm:pfacenotsplit}).  

We start by extending the description of the set  $\Hom_R(F^e_*R,R)$ from the normal (see for example \cite{PayneFrobeniusSplitToric}) to the seminormal case.  Here $F_*^eR$ denotes the module $R$ with the action given by the Frobenius $F^e$. Then $F_*^eR$ can be identified with $R^{\frac{1}{p^e}}$, see \ref{sec:FeRdescription}. The basic building blocks of this set are given as follows. For a fractional lattice point $a\in \frac{1}{p^e}M$, we let 
\begin{equation}\label{eq:pi} \pi_a \colon F_*^ek[M] \to k[M],\  x^u \mapsto \left\{\begin{array}{ll}
x^{a+u} &\textrm{ if } a+u \textrm{ is in } M, \\
0 &\textrm{ otherwise.}\end{array}\right. 
\end{equation}
\begin{theorem}
 \label{thm:maps}
  Let $R$ be a seminormal monoid algebra. Suppose $e$ is so that $p^e$ strictly  bounds the $p$-power torsion of $M \cap \< D \> / M_D$ for any face $D \prec C$.  We have that $\pi_a \in \Hom_R(F^e_*R,R)$ if and only if the following conditions hold.
 \begin{enumerate}[label=(\arabic*), ref=(\arabic*)]
    \item  $\pi_a\in \Hom_R(F^e_*\overline{R},\overline{R})$, i.e.,  for the extremal rays $\rho \prec \sigma$ with primitive vectors  $v_{\rho}$, we have $\< a, v_{\rho}\> > -1$. \label{thmitem:-1}
    \item For all RUFs $D$ that are not properly contained in any $p$-face, we have: 
        \begin{enumerate}[label=(\alph*), ref=\theenumi(\alph*)]
            \item  If $a \in \< D \>$, then 
                \begin{enumerate}[label=(\roman*), ref=\theenumii{}(\roman*)]
                    \item  if $D$ is not a $p$-face, we have $a \in \frac{1}{p^e} M_D$; \label{thmitem:psat}
                    \item if $D$ is a $p$-face,  we have $a\notin \frac{1}{p^e} \widetilde{M}_D$, where \\ $\widetilde M_D = \{ x \in M \cap \< D \> \, | \, mx \in M_D \mbox{ for some } m \in \Z \mbox{ with } \gcd(m,p) = 1 \}$.\label{thmitem:punsat}
                \end{enumerate}
            \item If $a \not \in \< D \>$, then there exists $v\in D^*$ with $\<a , v\> > 0$.\label{thmitem:notinD}
        \end{enumerate}
\end{enumerate}
Moreover,  for any $p$-face $D$, the image of any $\pi_a\in \Hom_R(F^e_*R,R)$ does not intersect $k[D]$. 
 \end{theorem}

This theorem allows us to compute the $F$-splitting dimension, the splitting prime and the $F$-splitting ratio of a seminormal monoid algebra. 
The $F$-splitting ratio is a measure of singularities generalising the $F$-signature to the non-normal case (see Section \ref{sec:Fsplittingratiodef}). For seminormal monoid rings, it has been independently calculated by De Stefani, Monta\~no, and N{\'u}\~nez-Betancourt \cite[Proposition 4.6]{StefaniMontanoNunezBetancourt}. 
More precisely, they compute the regularity of a seminormal monoid ring and they introduce a new notion called a``pure monoid dimension" and prove that it is bounded above by the depth. The pure monoid dimension measures how far $R$ is from being normal. They introduce and compute a limit called the pure monoid ratio which coincides with the $F$-splitting ratio when $k$ is of prime characteristic. We restate the statement for the $F$-splitting ratio here, as we are using slightly different language.

Let $\sigma$ be the dual cone to $C$  in $N=\mathrm{Hom}(M,\Z)\otimes \bR$ and for $\rho$ an extremal ray in $\sigma$ let $v_{\rho}$ be the primitive generator of $\rho$.
To the normalisation of $S$, $\overline{S}$ is associated a rational polytope 
\[P_{\overline{S}} = \{ x \in M_{\R}  \mid
0 \leq \langle x,v_{\rho} \rangle < 1 \mbox{ for all  $\rho$}\}.\]

 For a normal monoid ring, the lattice volume of this polytope computes the $F$-signature \cite[Theorem 5.1]{WatanabeYoshidaMinimalRelative} (\textit{cf.} \cite{Korff}). For seminormal monoids the formula has to be modified as follows. 

\begin{theorem}\label{thm:Fsplittinratio}
 Let $S$ be a finitely generated pointed seminormal submonoid of a lattice  $M\cong \Z^n$. Assume that $R=k[S]$ is $F$-split.  Let 
 \[ D_S = \bigcap_{D\ \mathrm{RUF}} D\]
 be the intersection of all RUFs of $S$. 
The $F$-splitting ratio of $R$ equals the lattice volume (with respect to $M_{D_S}$) of \[P_S:= P_{\overline{S}} \cap D_S.\]
Moreover, the $F$-splitting numbers are given by 
\[
a_e(R) = \# \left( \frac{1}{p^e}S \cap P_S\right).
\]
Finally, the splitting dimension of $S$ is $\delta = \dim(D_S)$
and the splitting prime is the ideal generated by the monomials in the complement of $D_S$. 
\end{theorem}

See Section \ref{sec:Fsplittingratio} for our proofs. Note that $D_S$ is the F-pure face of \cite{StefaniMontanoNunezBetancourt}. We also show that for a seminormal monoid, the $F$-splitting ratio is bounded by the $F$-splitting ratio of the normalisation, see Proposition \ref{lem:splittingratioSignature}. 

Cartier fixed ideals, especially the test ideal, play an important role in singularity theory. 

In Section \ref{sec:testideals} we use Theorem \ref{thm:maps} to calculate the Cartier fixed ideals. Let $\sC$ denote the Cartier algebra (see Section \ref{sec:Cartier} for definitions). We  can describe the $\sC$-fixed ideals of $R=k[S]$ in Theorem  \ref{thm.fixedideals}.

\begin{corollary}\label{cor:Cfixed} Let $R=k[S]$ be a seminormal affine monoid algebra. 
\begin{enumerate}
        \item The $\sC$-fixed ideals are radical and closed under intersection. 
        \item  The test ideal is given by 
        \[ \tau (R) = \< x^u \mid u \in S \text{ and u is not on any proper RUF } \>\]
        \item The non-$F$-pure ideal is given by
  \[
  \sigma(R) = \langle x^u \mid u \in S \mbox{ and $u$ is not on any $p$-RUF} \rangle.
  \]
    \end{enumerate}
\end{corollary}

The following example is a first illustration of our results. 

\begin{example}\label{ex:running} The Whitney umbrella. 
Let $S\subset \bZ^2$ be the submonoid generated by $\{(2,0),(0,1),(1,1)\}$, see Figure \ref{fig:ex1}. We have $k[S]=k[x^2,y,xy]=k[u,v,w]/(uv^2-w^2)$. Then $C$ is the rational polyhedral cone in $\bR^2$ spanned by $(1,0)$ and $(0,1)$, and it has 4 faces, $C$, $\<(1,0)\>$, $\<(0,1)\>$, and $\{(0,0)\}$ of which $C$ and $\<(1,0)\>$ are RUFs. For $p\neq 2$, $k[S]$ is $F$-split, so fix $p\neq 2$. We have $D_S = \<(1,0)\>$, the splitting dimension is $s=1$, and the splitting prime is the ideal generated by $x^{(0,1)}$ and $x^{(1,1)}$ in $k[S]$. $P_{\overline{S}}=\mathrm{conv} \<(0,0),(1,0),(0,1),(1,1)\>$ and $P_S = \mathrm{conv} \<(0,0),(1,0)\>$. As $M_{D_S}=2\bZ$, this implies that that the $F$-splitting ratio of $k[S]$ is $1/2$. The test   ideal agrees with the splitting prime and since $R$ is split, the non-$F$-pure ideal is the maximal ideal. 
    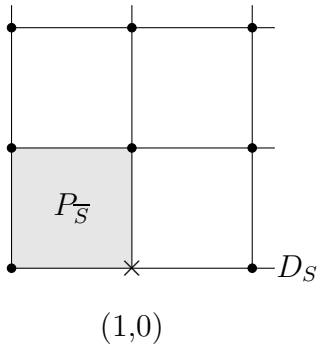
\begin{figure}[ht]\label{fig:ex1}
\begin{center}
\begin{tikzpicture}[scale = .4]
\fill[mygray] (0,0) rectangle (4,4);
\draw (2,2) node {$P_{\overline{S}}$};
 /
\draw (9.5,0) node {$D_S$};
\draw (0,8.75) -- (0,0) -- (8.75,0);
\draw (8,8.75) -- (8,0);
\draw (4,8.75) -- (4,0);
\draw (0,8) -- (8.75,8);
\draw (0,4) -- (8.75,4);
	\draw (4,0) node {$\times$};
	\foreach \i in {0,2,4}
	\foreach \j in {2,4}
		\filldraw [black] (2*\i,2*\j) circle (4pt);
	\filldraw [black] (0,0) circle (4pt);
	\filldraw [black] (8,0) circle (4pt);
\draw (4,-2) node {(1,0)};
\end{tikzpicture} 
\end{center}

\begin{center}
\caption{The monoid $S$ Example~\ref{ex:running}.}
\end{center}
\end{figure}
\end{example}

\noindent
  \textbf{Acknowledgements:}
    We would like to thank Johannes Hofscheier.

\section{Preliminaries}
\label{sec:Preliminaries}
\subsection{Seminormality and weak normality}\label{sec:basics}
Throughout this article, we work over a perfect field $k$ of characteristic $p > 0$.  When $e \in \N$ is understood, we write $q = p^{e}$ to denote the corresponding power of $p$. Recall that a Noetherian domain $R$ with fraction field $K$ is said to be \emph{seminormal} \cite{TraversoPicardGroup, SwanSeminormality} if whenever $x \in K$ satisfies $x^2,x^3 \in R$ one has $x \in R$. Equivalently, for any $x \in K$, one has $x \in R$ provided $x^N \in R$ for all $N \gg 0$. Similarly, $R$ is said to be \emph{weakly normal} \cite{YanagiharaWeaklyNormal} if whenever $x \in K$ satisfies $x^p \in R$ one has $x \in R$, or equivalently $x \in R$ provided $x^q \in R$ for every $q$. One sees immediately that weak normality implies seminormality. 

The seminormalization of a ring can also be characterized geometrically as the largest partial normalization that remains a bijection on points, and seminormality plays a prominent role in the study of (notably semi-log-canonical) singularities in complex birational algebraic geometry (see \cite{KollarSingsofMMP}).

Let $R$ be a ring with positive characteristic $p > 0$ and $F$ denote the Frobenius morphism. Let $F^e_*R$ be the module $R$ with $R$ action via $r_1\cdot r_2= r_1^qr_2$. It is often useful to visualize $F^e_*R$ as $R^{\frac{1}{q}}$, and then Frobenius is the inclusion $F^e: R\to F^e_*R$. A splitting is an $R$-module homomorphism $\pi: F_*R \to R$ such that $\pi \circ F = \mathrm{id}_R$ and  $R$ is called $F$-split if a splitting exists. A $F$-split ring is necessarily reduced, and moreover must be weakly normal \cite{HochsterRoberts}, so also seminormal.

\subsection{Definition of \texorpdfstring{$F$}{F}-splitting ratio used in Section \ref{sec:Fsplittingratio}}\label{sec:Fsplittingratiodef}
An $F$-splitting necessarily induces a rank one free direct summand of $F_*^e R$, and similarly for the iterates of Frobenius as well.  By a fundamental result of Kunz \cite[Lemma 51.17.6]{stacks-project} recall that a positive characteristic ring is regular if and only if the Frobenius endomorphism is flat.  In the case of $F^e_* R$ for an $F$-finite local ring $R$, this means that $R$ is regular if and only if $F^e_* R$ is a finite rank free module. As such, it is natural to analyze singularities by measuring how close $F^e_* R$ is to being free.

\begin{definition}
If $R$ is a ring of characteristic $p > 0$ and $e > 0$, the \emph{$e$th
$F$-splitting number of $R$} is the maximal rank $a_{e}(R)$ of a free
$R$-module direct summand of $F^{e}_{*}R$.  $R$ is said to be
$F$-split if $a_{e}(R) > 0$ for some (equivalently all) $e > 0$.
\end{definition}

\begin{theorem}
  If $R$ is an $F$-split local or standard graded ring with perfect residue field,
  there exists a positive integer $\delta$ such that the limit
\[
r_{F}(R) := \lim_{e \to \infty} \frac{a_{e}(R)}{p^{e\delta}}
\]
exists and $0 < r_{F}(R) < \infty$.  In this case, $\delta$ is said to
be the \emph{$F$-splitting dimension of $R$} and $r_{F}(R)$ is
called the \emph{$F$-splitting ratio of $R$}. The $F$-splitting dimension is alternatively characterized as the dimension of the prime ideal
\[
P = \langle f \in R \mid  \phi(F^e_*f) \in \mathfrak{m} \mbox{ for all $e > 0$ and all } \phi \in \Hom_R(F^e_*R, R) \rangle
\]
which is called the \emph{$F$-splitting prime}. When $P = 0$ or equivalently $\delta = \dim(R)$, $r_F(R)$ is denoted $s(R)$ and called the \emph{$F$-signature of $R$}.
\end{theorem}

The study of the growth rate of the ranks of the free direct summands of $F^e_* R$ originated in \cite{SmithVanDenBerghSimplicityOfDiff}, and the term $F$-signature was coined in \cite{HunekeLeuschkeTwoTheoremsAboutMaximal} where the existence of the limit in the local Gorenstein case was shown.
The $F$-splitting ratio, dimension, and prime were first defined in \cite{AberbachEnescu} who also conjectured the precise relationship formulated above.  The full version of the above result was shown in \cite{TuckerFsignExists,BlickleSchwedeTuckerFSigPairs1} (see also \cite{PolstraTucker} for a number of simplified arguments).  More recently, a global analysis of the $F$-splitting numbers and $F$-signature has also been carried out in \cite{DeStefaniPolstraYaoGlobalFsig}, though we have opted for the simplified version above in our exposition as it will suffice in the setting of monoid algebras.

\subsection{Cartier algebra and test ideals used in Section \ref{sec:testideals}.}\label{sec:Cartier}
Let $R$ be a commutative Noetherian ring of characteristic $p > 0$.  Recall that the \emph{Cartier algebra of $R$}, at times in the literature also called the \emph{total algebra of $p^{-e}$-linear maps on $R$}, is the $\mathbb{N}$-graded ring $\sC = \bigoplus_{e \geq 0} \sC_e$ where $\sC_e = \Hom_R(F^e_*R,R)$ is the set of all potential splittings of the $e$-iterated Frobenius.  Multiplication of homogenous elements in the Cartier algebra is given by function composition: if $\phi_1 \in \sC_{e_1}$ and $\phi_2 \in \sC_{e_2}$, then 
\[
\phi_1 \cdot \phi_2 = \phi_1 \circ F^{e_1}_* \phi_2 \in \Hom_R(F^{e_1 + e_2}_{*}R, R) = \sC_{e_1 + e_2}.
\]
See \cite{BlickleTestIdealsViaAlgebras,SchwedeTestIdealsInNonQGor} for further details.

\begin{definition}
Let $\sC = \bigoplus_{e \geq 0} \Hom_R(F^e_*R,R)$ denote the Cartier algebra of $R$.
An ideal $I \subseteq R$ is said to be \emph{$\sC$-fixed} if
\begin{equation}
\label{eq.fixed}
I = \sum_{e >0} \sum_{\phi\in \Hom_R(F_*^eR,R)} \phi (F^e_* I).
\end{equation}
The \emph{test ideal} $\tau(R)$ is the smallest non-zero $\sC$-fixed ideal of $R$, and the \emph{non-$F$-pure ideal} $\sigma(R)$ is the largest  $\sC$-fixed ideal of $R$.
\end{definition}
\begin{remark}
An ideal $I$ is called $\mathcal{C}$-compatible if 
  $I \subseteq \sum_{e >0} \sum_{\phi\in \Hom_R(F_*^eR,R)} \phi (F^e_* I)$. If $R$ is $F$-split, then a $\mathcal{C}$-compatible ideal is $\mathcal{C}$-fixed. 
\end{remark}

The following is well-known (see \cite[Proposition 2.3]{HsiaoSchwedeZhangCartieronToric}).

\begin{proposition}
\mbox{ }
\begin{enumerate}
 \item
 The sum of $\sC$-fixed ideals is also $\sC$-fixed.
 \item
 The non-$F$-pure ideal $\sigma(R)$ exists, and $\sigma(R)=R$ if and only if $R$ is $F$-split.
 \item
 If $R$ is $F$-split, then the intersection of $\sC$-fixed ideals is $\sC$-fixed.
 \item
 If $R$ is $F$-split, then any $\sC$-fixed ideal is a radical ideal and all of the minimal primes of $\sC$-fixed ideals are also $\sC$-fixed.
 \item
 If $c \in R$ is such that $R_c$ is regular, there is some $N > 0$ so that $c^N$ belongs to every non-zero $\sC$-fixed ideal.
 \item
 The test ideal $\tau(R)$ exists.
 \end{enumerate}
\end{proposition}

\begin{remark}
Note that $R$ is strongly $F$-regular if and only if $\tau(R) = R$. 
 When $R = k[S]$ is an affine monoid algebra, strongly $F$-regular is equivalent to the condition that $R$ is normal.

\end{remark}

\section{Seminormal Monoid Algebras}\label{sec:Seminormal monoids}

Seminormal monoid rings have been studied by a number of researchers, including notably Bruns, Li, and R{\" o}mer in \cite{BrunsLiRoemer} and by Li in her thesis \cite{Lithesis}. 

An \emph{affine monoid} $S$ is a finitely generated commutative 
monoid which can be embedded into a lattice $\Z^n$ for some $n\in \N$. We always use $+$ for the monoid operation and set $M = \Z S$. 

To $S$ is associated the rational polyhedral 
convex cone $C(S) = \R_+S$ in $M_{\R} = M\otimes_{\Z}\R$.
 
The saturation $\overline{S} = C \cap M$ of $S$ is also
a finitely generated monoid. 
We assume that $C(S)$ is a pointed cone and say equivalently that $S$ is a positive monoid. Note that if $C$ is not pointed, we can reduce to the pointed case by quotienting by the largest subspace contained in $C$. 
If $R = k[S]$ is the associated monoid algebra
over a field $k$,  then $\overline{R} = k[\overline{S}]$ is
the normalization of $R$.

For a face $D \prec C$, let $M_{D} = \Z (S \cap D)$ denote the sublattice of $M$ generated by $S\cap D$. 
Note that if $D'\prec D$, we have $M_{D'}\subset M_{D}$. Moreover, if $\langle D \rangle = \R D$ denotes the linear span of $D$, then $\rank M_D = \rank M \cap \< D \> = \dim \< D \> $.  In particular, for each face $D \prec C$, we have that $ M \cap \< D \> / M_D$ is a finite abelian group. For a face $D$, we denote its interior by $\mathrm{int}(D)$.

\begin{definition}\label{def:seminormal}
   We say that $S$ is \emph{seminormal} if for all $D\prec C$ we have 
\[ S \cap \mathrm{int}(D) = M_{D} \cap \mathrm{int}(D).\]
\end{definition}

\noindent
Observe that $M_C = M$ as we have assumed $M = \Z S$.  Combining \cite[Theorem 3.8]{ReidRoberts} (c.f., \cite[Proposition  5.32]{HochsterRoberts}) with \cite[Theorem 4.3]{ReidRoberts}, we have the following.
\begin{proposition}
The affine monoid $S$ is seminormal if and only if the associated monoid algebra $k[S]$ is seminormal (for any field $k$). 
\end{proposition}

Throughout this article, unless explicitly stated otherwise, we will assume that the monoid $S$ is seminormal.

\begin{definition}\label{def:NNS}
We say that a proper face  $D \precneqq C$  is \emph{relatively saturated} if  
\[M_{D} = \bigcap _{D\precneqq D'} M_{D'} \cap \<D \>,\] where the 
intersection is over all faces $D'$ strictly containing $D$. If $D$ is 
not relatively saturated, we say that $D$ is a \emph{relatively unsaturated face}  (RUF). 
We also use the convention that $C$ is a RUF and call a RUF that is not $C$ a \emph{proper RUF}. 
\end{definition}

The following proposition gives a convenient way to check which elements in the normalization of $S$ are in $S$.

\begin{proposition}\label{prop:inS}
	Let $S$ be a seminormal monoid. 
	\begin{enumerate}
		\item \label{proppart1} If $u\in \overline{S}$, then $u\in S$ if and 
		only if the following condition is satisfied: for all RUF's $D$ with $u \in D$, we have $u \in M_D$.
		\item \label{proppart2} Moreover, for any face $D \prec C$, we have
		\begin{equation*}
			M_D = \bigcap_{\substack{D \prec D'  \\ D' \mathrm{RUF}}} M_{D'} \cap \< D\>.
		\end{equation*}  
	\end{enumerate}
\end{proposition}

\begin{proof}
	We first show that
	\begin{equation}
	\label{eq:claiminproof} \tag{$*$} \textstyle
		\bigcap _{D\precneqq \widetilde{D}} M_{\widetilde{D}}\cap \langle D\rangle = 
	 \bigcap _{\substack{D\precneqq D' \\ D'\ \mathrm{RUF}}} M_{D'}\cap \langle D\rangle.
	\end{equation} 
	for all proper faces $D\precneqq C$.
	If $\mathrm{codim}(D) = 1$,
	then $C$ is the only strictly larger face and is relatively unsaturated, so both sides are just $M \cap \langle D \rangle$. 
	We now proceed by induction and assume that $\mathrm{codim}(D) > 1$. 
	Note that
	\begin{equation}
		\label{eq:tempstep} \tag{$**$}
		\textstyle
			 \bigcap _{D\precneqq \widetilde{D}} M_{\tilde{D}}\cap \langle D\rangle =
			  \left(\bigcap _{\substack{D\precneqq D' \\ D'\ \mathrm{RUF}}} M_{D'}\cap \langle D\rangle \right) \bigcap \left(\bigcap _{\substack{D\precneqq \widetilde{D} \precneqq C \\ \widetilde{D} \mbox{ \footnotesize rel. sat.}}} M_{\widetilde{D}}\cap \langle D\rangle \right).
	\end{equation}
	For each relatively saturated face $D\precneqq \widetilde{D} \precneqq C$, we have by substituting in from the definition of relatively saturated and using the induction hypothesis 
	\begin{equation*}
		\begin{array}{rcl}
		M_{\widetilde{D}} \cap \langle D \rangle &=& \left( \bigcap _{\widetilde{D}\precneqq D''} M_{D''} \cap \<\widetilde{D} \> \right) \cap \langle D \rangle = \left( \bigcap _{\substack{\widetilde{D}\precneqq D'' \\ D''\ \mathrm{RUF}}} M_{D''} \cap \<\widetilde{D} \> \right) \cap \langle D \rangle \\
		&=& \bigcap _{\substack{\widetilde{D}\precneqq D'' \\ D''\ \mathrm{RUF}}} M_{D''} \cap \<D \>  \subseteq \bigcap _{\substack{D\precneqq D' \\ D'\ \mathrm{RUF}}} M_{D'}\cap \langle D\rangle.
		\end{array}
	\end{equation*}
	Thus, the intersections on the right side in \eqref{eq:tempstep} are redundant, and so \eqref{eq:claiminproof} follows.
	
	To prove part (\ref{proppart1}) of the proposition, 
	suppose first $u\in S$. Then if $D$ is any face containing $u$, we have $u\in S\cap D \subset M_D$ by definition of $M_D$. 
	For the reverse implication, now assume that $u$ is contained in $M_D$ for all RUFs containing $u$. Note that $S=\bigcup _{D\preceq C}(M_D \cap \mathrm{int}(D))$ by Definition \ref{def:seminormal}. Let ${D}$ be the smallest face containing $u$, so that $u\in \mathrm{int}(D)$. We have to show that $u\in M_{{D}}$. If ${D}$ is a RUF, then $u\in M_{{D}}$ by assumption. Otherwise, assume that ${D}$ is a relatively saturated face, and so by definition $M_{{D}} = \bigcap _{{D} \precneqq \widetilde{D}} M_{\widetilde{D}} \cap \langle D \rangle$. The statement then follows from the assumption and \eqref{eq:claiminproof}.

	Finally, for the equation in part (\ref{proppart2}) of the proposition, note that if $D$ is a RUF both sides agree as $D$ necessarily appears as part of the intersection on the right. Else, if $D$ is relatively saturated it is necessarily proper and the statement follows from \eqref{eq:claiminproof} and the definition of relatively saturated.
	\end{proof}

For a monoid algebra, weakly normal is in fact equivalent to being F-split, and can be characterized using the following variation on the relatively unsaturated condition.

\begin{definition}
\label{def:pface}
We say that a face $D \preceq C$ is a \emph{$p$-face} if $p$ divides $[M \cap \<D \> : M_D]$.  A $p$-face that is also relatively unsaturated will be abbreviated as a pRUF.  If $D \preceq C$ is a pRUF that is not properly contained in any other pRUF, we say that $D$ is a maximal pRUF.
\end{definition}

\begin{theorem}\label{thm:pfacenotsplit}
\cite[Proposition 6.2]{BrunsLiRoemer}
For a seminormal monoid $S$, the following are equivalent:
\begin{enumerate}
    \item  $R = k[S]$ is $F$-split; 
    \item $R = k[S]$ is weakly normal; 
    \item $C(S)$ has no $p$-faces. 
\end{enumerate}
\end{theorem}

The following lemma shows that maximal $p$-faces must be pRUFs. 

\begin{lemma}\label{lem:pfaces}
Suppose that $S$ is seminormal, and $D \precneqq C$.
\begin{enumerate}
 \item
 If $D$ is a $p$-face, then $D$ is contained in some maximal pRUF.
\item
$D$ is a maximal $p$-face, i.e., any stricly larger face $D \precneqq D'$ is not a $p$-face, if and only if  $D $ is a maximal pRUF.
\item
$D$ is not contained in any maximal pRUF if and only if all strictly larger faces $D \prec D'$ are not $p$-faces.
\end{enumerate}
 \end{lemma}

\begin{proof}
\textit{(1)} If $D$ is a RUF and hence pRUF, then it is clearly contained in a maximal pRUF.  Else, if $D$ is relatively saturated, then we have that $M_D = \bigcap_{\substack{D \prec D'  \\ D' \mathrm{RUF}}} M_{D'} \cap \< D\>$.  By assumption, there exists some $x \in M \cap \<D \> \setminus M_D$ with $px \in M_D$.  Since $x \not\in M_D$, there exists some $D \prec D'$ with $D'$ a RUF and $x \not\in M_{D'}$.  But then $x \in M \cap \< D' \> \setminus M_{D'}$ and $px \in M_{D'}$.  Thus, we must have $p | [M \cap \< D' \> : M_{D'}]$ so that $D'$ is a pRUF.  As before, we conclude $D$ is contained in some maximal pRUF.

\textit{(2)} 
If $D$ is a maximal $p$-face, it must be contained in a $p$-RUF by (1), since it is a $p$-face. By maximality, it must be equal to it, so $D$ is a pRUF. Conversely suppose $D$ is a maximal pRUF contained in a $p$-face $\widetilde{D}$. By (1) $\widetilde{D}$ is contained in a maximal pRUF $D'$. Since $D$ is a maximal pRUF, we have $D=\widetilde{D}=D'$, so $D$ is a maximal $p$-face.

\textit{(3)} If $D$ is contained in a pRUF, then that pRUF is necessarily a $p$-face.  Conversely, if $D$ is contained in a $p$-face, then it is contained in a pRUF by \textit{(1)}.
\end{proof}

\section{The set of homomorphisms from \texorpdfstring{$F_*R$}{F*R} to  \texorpdfstring{$R$}{R}.}

\subsection{The action of Frobenius on \texorpdfstring{$R$}{R} and a description of \texorpdfstring{$F_*^eR$}{F*R}}\label{sec:FeRdescription}
As $S$ is
finitely generated,  it follows that $R = k[S]$ is a finitely generated $k$
algebra.  As such, when $k$ has characteristic $p > 0$, the
Frobenius or $p$th
power map $F : R \to R$ given by $r \mapsto r^{p}$ is a finite ring
endomorphism.  In the case of a monoid algebra, it is often convenient to view Frobenius as the
composite of two interchangeable steps, namely dilation by $p$ of the
underlying lattice $M$ and the application of Frobenius to the
underlying field $k$.  When $k$ is perfect, as we shall largely assume
throughout, this latter process is in fact an isomorphism. 
Similar remarks apply for the
iterates $F^{e} : R \to R$ of Frobenius for $e > 0$.

Given an $R$ module $\sM$, one may consider the modules $F^{e}_{*}\sM$
for $e > 0$ given by restriction of scalars from $F^{e}$.  In other
words, $F^{e}_{*}\sM$ agrees with $\sM$ as an abelian group but differs
as an $R$-module in that $R$ acts on the former via $q$-th powers, where we recall $q=p^e$.
Note that $F^{e}_{*}(\blank)$ is an exact functor on the category
of $R$-modules that commutes with arbitrary localizations.  Furthermore $F^{e}_{*}R$ is naturally identified with
the ring
$R^{1/q}$ of $q$-th roots of $R$.  From this perspective, the
$e$-iterated Frobenius $F^{e}$ takes on the guise of the inclusion $R
\subseteq R^{1/q}$.
Similarly, if $\sM$ is any $R$-submodule of
$K = \Frac(R)$, we have $F^{e}_{*}\sM \simeq \sM^{1/q}$.
When furthermore $\sM$ is an $M$-graded $R$-submodule of $k[M] \subseteq K$, we
may write $\sM = \bigoplus_{u \in U} k \cdot x^{u}$ for some
$S$-invariant subset $U \subseteq M$.  In this case, we may further
identify
\[
F^{e}_{*}\sM \simeq \sM^{1/q} = \bigoplus_{u \in \frac{1}{q}U} k \cdot
x^{u} 
\]
inside of $k[\frac{1}{q}M] \subseteq K^{1/q} 
\simeq F^{e}_{*}K$.  In particular, this
discussion applies to both $R$ and $\overline{R}$, where we have
$R^{1/q} = k[\frac{1}{q}S]$ and $\overline{R}^{1/q} = k[\frac{1}{q}\overline{S}]$. Though our main purpose is to study $F_*^eR$, we shall first need a description of $\mathrm{Hom}(F_*^eR, R)$ below. 

\subsection{Notation and review of the normal case.}

For normal affine monoid rings, the description of  $\mathrm{Hom}(F_*^eR, R)$ has appeared for example in \cite{PayneFrobeniusSplitToric}, . We will review  this description  and extend it to
seminormal affine monoid rings. Let $T = \mathrm{Spec}(k[M])$ be the open dense torus contained in $\mathrm{Spec}(R)$. Recall the definition of $\pi_a$ from equation \ref{eq:pi}. 
 
The set of these maps forms
a basis for $\mathrm{Hom}_{k[M]}(F_*^ek[M], k[M])$ as a vector space over $k$
by \cite[Lemma 4.1]{PayneFrobeniusSplitToric}\footnote{Suitably interpreted, the same statement remains true over $\Z$ as well.} and by the following Lemma they form the building blocks of 
$\mathrm{Hom}(F_*^eR, R)$.

Let $\sigma \in N_{\R}$ be the dual cone to $C$, and for a ray $\rho$ of $\sigma$,  let $v_{\rho}$ denote the primitive generator of $\rho$. Let $D_{\rho}=C\cap \rho^{\perp}.$ Recall that for a face $D\prec C,$ \[D^{\vee}=\{u\in \sigma \mid \< u,v\> \geq 0\ \forall\ v\in D\}\] is the dual face of $D$. 
\begin{definition}\label{def:D*}
For a face $D\prec C$ we let 
 $D^*=\sigma\cap D^{\perp} \prec \sigma$.
 \end{definition}
 Note that since $C$ is full-dimensional, $\sigma$ is pointed. This in turn implies that for all rays $\rho$ in $\sigma$, we have $(D_{\rho})^* = \rho$. 
We then have 
\begin{equation}\label{eq:relint}\relint(D)  = 
    \{u \in D \mid \<u,v_{\rho}\>>0 \textrm{ for all primitive vectors }v_{\rho}\in D^{\vee}\smallsetminus D^*\}.\end{equation}

\begin{lemma}\label{lem:sam}
The set of $\pi_a$ for all $a \in \frac{1}{q}M$ such that $(a+\frac{1}{q}S)\cap M \subset S$ forms a $k$-basis for $\Hom_R(F_*R,R)$. When
$S$ is normal, this set is given by $\{a \mid \langle a, v_{\rho}\rangle > -1\textrm{ for all rays } \rho \in \sigma\}$.
\end{lemma}

\begin{proof}
The statement is well-known (\textit{cf.} \cite[discussion in beginning of Section 4]{PayneFrobeniusSplitToric}).  For completeness, and because it is essential to all of our arguments below, we sketch an argument.  For the first assertion, note that any $\phi \in \Hom_R(F^e_*R,R)$ determines an extension to $\Hom_{k[M]}(F^e_*k[M],k[M])$, also written $\phi$ and having the property that $\phi(F^e_*R)\subseteq R$.  By \cite[Lemma 4.1]{PayneFrobeniusSplitToric}, we can then write the extension in terms of the $k$-basis described in equation \eqref{eq:pi} $\phi = \lambda_1 \pi_{a_1} + \cdots + \lambda_s \pi_{a_s}$ where $a_1, \ldots, a_s \in \frac{1}{q}M$ are pairwise distinct and $\lambda_1, \ldots, \lambda_s \in k$. For any $u \in \frac{1}{q}S$, as $\phi(x^u) = \sum_i \lambda_i \pi_{a_i}(x^u) \in R$ and $a_i + u \neq a_j + u$ for $i \neq j$ so that the $\pi_{a_i}(x^u)$ are distinct monomials when non-zero, it follows we must have $\pi_{a_i}(x^u) \in R$ for all $i$.  Varying $u \in \frac{1}{q}S$, we see $(a_i+\frac{1}{q}S)\cap M \subset S$ for all $i$ as desired.

For the second assertion, first note that it is easy to check that $a$ with $\<a, v_{\rho}\>> -1 \ \forall \ \rho \in \sigma$ satisfies $(a+\frac{1}{q}S)\cap M\subset S$. Now suppose we have some $a \in \frac{1}{q}M$ with $\langle a, v_\rho \rangle \leq -1$ for some $v_\rho$.  From the above description, it suffices to show that $\pi_a$ does not satisfy $(a + \frac{1}{q}S) \cap M \subseteq S$. 
We will need the following claim: 

\begin{tabular}{cp{13cm}}
($\dag$)&
Let $a\in \frac{1}{q}M$ and fix $\rho$.  Then there is $w\in M$ such that $\<w,v_{\rho}\> = -1$ and $\langle w, v_{\rho'} \> > \langle a, v_{\rho'} \rangle $ for all $\rho' \neq \rho$. 
\end{tabular}

Indeed, let $\widetilde{w}\in M$ such that $\<\widetilde{w},v_{\rho}\> = -1$. Let $w'\in \mathrm{relint}(D_{\rho}).$ By \eqref{eq:relint} $\<w',v_{\rho'}\>>0$ for all $ \rho'\neq \rho$ and $\<w',v_{\rho}\> =0$. So we can find a multiple $m$ such that $w:=\widetilde{w}+mw'$ satisfies the desired properties. 

 Now let $w$ be as in the claim and consider $u:= w - a \in \frac{1}{q}M$.  We have that $\langle u, v_\rho \rangle = -1 - \langle a, v_\rho \rangle \geq 0$ and $\langle u, v_{\rho'} \rangle \geq 0$ for $\rho \neq \rho'$ by construction, so it follows $u\in \frac{1}{q}S$. But $u+a =w \in \left( (a + \frac{1}{q}S) \cap M \right) \setminus S$ as desired.
\end{proof}

\subsection{The seminormal case}

The following Lemma is true in general \cite[Exercise 1.2.(4)]{BrionKumar05}. 
However, in our situation it can be seen directly 
using the combinatorial description of $R$. 

\begin{lemma}\label{lem:extensiontosaturation}
Let $R$ be a seminormal monoid algebra. Then every homomorphism $\pi_a \in 
\Hom(F_*^eR,R)$ extends uniquely to a homomorphism on the normalisation
$\pi_a \colon F_*^e\overline{R} \to \overline{R}$. 
\end{lemma}

\begin{proof}
First, since $C$ is full-dimensional, $k[T]$ is the ring of fractions of $R$. So a morphism 
$F_*^eR\to R$ extends uniquely to a morphism $F_*^ek[T] \to k[T]$. Thus it suffices 
to show that if $\pi_a\in \mathrm{Hom}(F_*^ek[T],k[T])$ restricts to a homomorphism $F_*^eR \to R$, then it restricts to a homomorphism $F_*^e\overline{R} \to \overline{R}$. 
In other words, by Lemma~\ref{lem:sam}, we want to show that if $\pi_a \in \Hom(F_*^eR,R)$, then $\<a, v_\rho\> > -1$ for all rays $\rho$ in $\sigma$.

Assume there is $\rho$ such that $\<a,v_{\rho}\> \leq 1$. We need to find $u\in \frac{1}{q}S$ such that $\pi_a(x^u) \notin R$, i.e., $u+a\in M\smallsetminus S$. Let $D_{\rho}\prec C$ be the facet corresponding to $\rho$. 
By ($\dagger$)  in the proof of Lemma~\ref{lem:sam}, we can find $w\in M$ such that 
$\langle w, v_\rho \rangle = -1$ and $\langle w, v_{\rho'} \rangle > \langle a, v_{\rho'} \rangle $ for all $\rho' \neq \rho$. 
Let $u=w-a$. Clearly $w=u+a\in M \smallsetminus S$. It remains to show $u\in \frac{1}{q}S$. If $\<a, v_{\rho}\> <-1$, then $u\in\frac{1}{q}M\cap \mathrm{int}(C)\subset\frac{1}{q} S$. If $\<a, v_{\rho}\> = -1$ and $D_{\rho}$ is relatively saturated, 
then $u\in \frac{1}{q}M \cap \mathrm{int}(D_{\rho}) = \frac{1}{q}M_{D_{\rho}}\cap 
\mathrm{int}(D_{\rho}) \subset \frac{1}{q}S$. 
If  $\<a, v_{\rho}\> = -1$ and $D_{\rho}$ is relatively unsaturated, let $\widetilde{w}\in M \cap \<D_{\rho}\>\smallsetminus M_{D_{\rho}}$ with $\<\widetilde{w}, v_{\rho'}\> > \< a, v_{\rho'}\> $ for all rays $\rho' \neq \rho$ in $\sigma$. Then take $u= \widetilde{w} -a \in \mathrm{int}(C) \cap \frac{1}{q}M \subset \frac{1}{q}S$.  
\end{proof}

 \begin{proof}[Proof of Theorem~\ref{thm:maps}]
 Assume $\pi_a \in \Hom(F_*^eR,R)$. 
 First note that \ref{thmitem:-1} follows from Lemma~\ref{lem:extensiontosaturation}.

    Let $D$ be a RUF that is not properly contained in a $p$-face. 
       To see \ref{thmitem:psat},
     by way of contradiction, assume that $\<D\>$ is not a $p$-face and we have $a\in \<D\>\smallsetminus \oq M_D$. 
    We claim that there is $w\in M\cap \<D\>\smallsetminus M_D$ such that $qw-qa\in M_D$. Indeed, since $p\nmid [M\cap\<D\> \colon M_{D}]$, multiplication by $p$ on $M\cap \<D\>/M_D$ is an isomorphism. So there exists $w\in M\cap \<D\>$ such that $qw+M_D = qa+M_D$. Since $qa\notin M_D$, we must have $w\notin M_D$, and the claim follows. 
  Now pick $w'\in M_D\cap \relint(D)$. There is $\ell\in \Z_{>0}$ such that $u:=-a+w+\ell w'\in \relint (D)$ and $w+\ell w'\in \relint (D)$. Then $u\in \oq M_D\cap \relint(D) \subset \oq S$, but 
 $u+a=w+\ell w'\in M\cap \relint(D) \smallsetminus M_D$, so 
 $(a+\oq S)\cap M\subsetneq S$, a contradiction to $\pi_a \in 
 \Hom_R(F^e_*R,R)$.
    
     We now show \ref{thmitem:punsat}.  By way of contradiction, assume that $\<D\>$ is a $p$-face and $a\in \<D\>\cap \oq \widetilde{M}_D$. We first show that the assumption that  $\pi_a\in \Hom(F_*^eR,R)$, then implies that $a$ must be in $\oq M_D$. We then show this is not possible if $D$ is a $p$-face. 
     
     By assumption, there is $m\in \Z$ with $\gcd(m,p)=1$, such that $mqa\in M_D$. Let $\alpha, \beta$ in $\Z$ such that $\alpha m+\beta q=1$, so $\alpha ma +\beta qa=a$. Let $w\in M_D\cap \relint (D)$ and $\ell \in \Z_{>0}$ such that $u:=\ell w-\alpha m a \in \relint(D)$ and $u+a =\ell w - \alpha ma +a \in \relint(D)$. Note that $u\in \oq M_D$, so $u\in \oq S$, and $u+a = \ell w+\beta q a\in M$. Since $\pi_a \in \Hom_R(F^e_*R,R)$, we have $(a+\oq S)\cap M \subset S$, 
     so we must have $u+a\in M_D$. But this implies that $a\in \oq M_D$. 
     
     Now, since $D$ is a $p$-face, $M\cap \<D\> /M_D$ has $q$-torsion, so there is $w\in M\cap \<D\>\smallsetminus M_D$ such that $qw\in M_D$. 
     Let $w'\in M_D\cap \relint(D)$ and pick $\ell$ sufficiently large such that $u=-a+w+\ell w'\in \relint(D)$ and $u+a = w + \ell w' \in \relint(D)$. But then $u\in \oq M_D\cap \relint(D)\subset \oq S$, but $u+a =w +\ell w'\in M \cap \relint (D) \smallsetminus M_D$, so $(a+\oq S)\cap M \subsetneq S$, a contradiction to $\pi_a \in \Hom_R(F^e_*R,R)$.

      For
      \ref{thmitem:notinD}, assume by way of contradiction that $a\notin \<D\>$ and $\<a,v\>\leq 0$ for all $v\in D^*$. Let $D'\prec C$ be the smallest face such that $\<D'\>$ contains $D$ and $a$. 
    We claim $a\in \oq M_{D'}\cap \<D'\>$. If $D'$ is an RUF, then by assumption we have $D'$ is not a $p$-face so (2)(a)(i) implies $a\in \oq M_{D'}$. If $D'$ is not an RUF, then for any  RUF $\widetilde{D}$  containing $D'$ we conclude in the same manner $a\in \oq M_{\widetilde{D}}\cap\<\widetilde{D}\>$ and since $M_{D'} = \cap_{D'\prec \widetilde{D}}M_{\widetilde{D}}\cap \<D\>$ by \ref{prop:inS} (2), it follows that $a\in \oq M_{D'}\cap\<D'\>$. 
      Since $D$ is a RUF, there is  $w\in M_{D'}\cap \<D\> \smallsetminus M_D$.  Let $w'\in M_D\cap \relint(D)$.
      We claim that there is $\ell \in \Z_{>0}$ large enough such that 
       $w+\ell w'\in \relint(D)$ and $u:=-a+w+\ell w' \in \relint(D')$.
       Then $u\in \oq M_{D'}\cap \relint(D')\subset \oq S$, but $a+u = w+\ell w'\in \relint(D)\cap M\smallsetminus M_D$. So $(a+\oq S)\cap M$ is not a subset of $S$, a contradiction to $\pi_a\in \mathrm{Hom}_R(F^e_*R,R)$. For the claim, note that for $\ell \in \Z_{>0}$ sufficiently large, we have $w+\ell w'\in \relint(D)$.  Now let $v_{\rho}\in 
       (D')^{\vee}\smallsetminus (D')^*=\left((D')^{\vee}\smallsetminus D^*\right)\cup \left(D^*\cap  (D')^{\vee} \smallsetminus (D')^*\right)$ be a primitive vector.  If  $v_{\rho}\in (D')^{\vee}\smallsetminus D^*$, we have  $\<w',v_{\rho}\>>0$ by \eqref{eq:relint} and we can choose $\ell$ sufficiently large so that $\<-a+w+\ell w',v_{\rho}\>>0$. 
    Now let $v_{\rho} \in D^*\cap (D')^{\vee}\smallsetminus (D')^*$. Note that by minimality of $D'$ and our assumption on $a$ we have $\<a,v_{\rho}\><0$ for all $v_{\rho}\in D^*\smallsetminus (D')^*$ . We also have 
     $\<w+\ell w',v_{\rho}\>=0$ for all $\ell$.   So $\<-a+w+\ell w', v_{\rho}\>>0$ for all $v_{\rho}\in D^*\cap (D')^{\vee}\smallsetminus (D')^*$. Now the claim follows from \eqref{eq:relint}. 
     
     To prove the reverse direction, we have to show that if $a$ satisfies the conditions of the theorem, and $u\in \oq S$ such that $\widetilde{u}:=a+u\in M$, then $\widetilde{u}\in S$. 
     
     We start by showing that $\widetilde{u}\in \overline{S}.$ Indeed, by \ref{thmitem:-1}, for any primitive generator $v_{\rho}$ of $\sigma$, we have $\<\widetilde{u}, v_{\rho}\> = \<a+u,v_{\rho}\> >-1.$ Since $\widetilde{u}\in M$, and $v_{\rho}\in N$, this implies $\<\widetilde{u},v_{\rho}\>\geq 0$, so $\widetilde{u}\in C$. 
     
     By Proposition~\ref{prop:inS}, it suffices to show that for every RUF $D$ containing $\widetilde{u}$, if $\widetilde{u}\in D$, then $\widetilde{u}\in M_D$.
    First note that 
     
    \begin{tabular}{cp{13cm}}
        ($\star$) & if $D$ is a RUF that is not properly contained in a  $p$-face, and $\widetilde{u}\in D$, then $a\in \<D\>$. 
\end{tabular}

Indeed, by \ref{thmitem:notinD}, if $a\notin \<D\>$, there is $v\in D^*$ such that $\<a,v\>>0$. Since $u\in C$, we have $\<u,v\>\geq 0$, so $\<\widetilde{u}, v\>>0$, a contradiction to $\widetilde{u}\in D$. 

We now show that $\widetilde{u}$ is not contained in any  $p$-face. Aiming for a contradiction, let $\widetilde{u}$ be contained in a  $p$-face  $D$,  and let $D'$ be a  maximal pRUF containing  $D$ as in Lemma \ref{lem:pfaces}. Since $\widetilde{u}\in D'$, we have $a\in D'$ by $(\star)$. Note that  $u=\widetilde{u}-a\in \<D'\> \cap \oq S$, so $u\in \oq M_{D'}\subset \oq \widetilde{M}_{D'}$ and $\widetilde{u}\in (a+\oq M_{D'})\cap M$. One can show that there exist positive integers $a_1, \ldots, a_m$ such that 
$M\cap \<D'\>/\widetilde{M}_{D'} \cong \bZ/(p^{a_1})\oplus \cdots \oplus \bZ/(p^{a_m})$, so we can conclude that 
$q\widetilde{u}\in \widetilde{M}_{D'}$, since $q$ bounds the $p$-power torsion of $M/\widetilde{M}_{D'}$.  But then  $a = \widetilde{u}-u\in \oq\widetilde{M}_{D'}$, a contradiction to the assumption \ref{thmitem:punsat} which says $a\notin \oq \widetilde{M}_{D'}$.

Now let $D$ be an RUF containing $\widetilde{u}$. 
By the previous paragraph, $D$ is not a $p$-face and not contained in any $p$-face. By $(\star)$ we have 
$a\in \<D\>$. By assumption \ref{thmitem:psat}, we have $a\in \oq M_D$.  Thus $\widetilde{u}\in \oq M_D\cap M\cap \<D\>$, in particular $q\widetilde{u}\in M_{D}.$ Since $D$ is not a $p$-face, $M\cap \<D\> /M_D$ has no $q$-torsion, so $\widetilde{u}\in M_D$. \end{proof}

Example \ref{Ex:needqbound} shows that in Theorem \ref{thm:maps} we cannot omit the assumption that $q$ bounds the $p$-power torsion without changing the statement. 

\begin{example}\label{Ex:needqbound}
Let $S$ be the submonoid generated by $(4,0),(0,1),(1,1),(2,1),(3,1)$, let $p=2$ and  $e=1$, so $q=p=2$. 
Let $D = \<(1,0)\>$. Then $M_{D} = 4\Z$, so $D$ is a $2$-face, and $\widetilde{M_D}=4\Z$. Note that  $a=(1,0)$ satisfies the conditions of the theorem. However, $(2,0) \in \frac{1}{2}S$, but $a+(2,0)\in M \smallsetminus M_D$. So $\pi_a\notin \Hom_R(F^e_*R,R)$.
\end{example}

\begin{question}\label{que:smallemaps}
Suppose $R=k[S]$ is not $F$-split. Describe the set $\Hom(F^e_*R,R)$ for small $e$. 
\end{question}

\subsection{Computation of the \texorpdfstring{$F$}{F}-splitting ratio and the \texorpdfstring{$F$}{F}-splitting prime.}\label{sec:Fsplittingratio}
We now proceed to 
compute  the $F$-splitting ratio of a seminormal monoid ring using the description of the homomorphisms from Theorem \ref{thm:maps}.

\begin{definition}
Let $P$ be a polytope and $M$ be a lattice. We denote by $\Vol_{M}(P)$ the volume of $P$ with respect to the lattice. 
\end{definition}
Note that if $M'\subset M$ is a sublattice of index $n$, then $\Vol_M(P) = n\Vol_{M'}P$.

\begin{proof}[Maps proof of Theorem \ref{thm:Fsplittinratio}] 
Recall that $a_e$ denotes the number of free summands in $F^e_*R$. 
We first claim that $a_e = | P_{\overline{S}} \cap D_S \cap \frac{1}{q}S|$. 
Let \[I(F^e_*R) = \left\langle x^u \in F^e_*R\mid \pi_a\left(x^{u}\right) \in \bm \mbox{ for all } a\in \frac{1}{q}M \mbox{ with } \pi_a \in \Hom_R(F^e_*R,R)\right\rangle.
    \] By an analogous statement to \cite[Proposition 4.5]{TuckerFsignExists} we have $\ell\left(F^e_*R/I(F^e_*R)\right)= a_e(R)$ for all $e$.
We have that $F^e_*R /I(F^e_*R)$ has as basis over $k$ the (images of the) set $x^u$ with $u \in \frac{1}{q}S$ and $\pi_a(x^u) \not\in \mathfrak{m}$ for some $\pi_a \in \Hom_R(F^e_*R,R)$. For such a basis element $u$ and $\pi_a \in \Hom_R(F^e_*R,R)$, we have $\pi_a(x^u)\notin \bm $ if and only if $a+u=0$, so $\pi_{-u}\in \Hom_R(F^e_*R,R)$, and $-u$ must satisfy the conditions of Theorem \ref{thm:maps}. Note that $R$ is $F$-split, so there are no $p$-faces and  Theorem \ref{thm:maps} applies for all $e$. Now $u\in \frac{1}{q}S$ implies $\langle u, v_{\rho}\rangle \geq 0$ for all $\rho$, and Theorem \ref{thm:maps} (1) for $-u$ implies that $\langle -u, v_{\rho}\rangle >-1$ for all $\rho$, so $u\in P_{\overline{S}}$. To see that $u\in D_S$, suppose there is a RUF $D$ such that $u\in P_{\overline{S}} \smallsetminus \langle D\rangle$. Then $-u \notin \<D\>$, so by Theorem \ref{thm:maps} (2) (b) there is $v\in D^*$ with $\<-u,v\>>0$, so there is some $\rho$ s.t. $\<u,v_{\rho}\><0$, a contradiction to $u\in \frac{1}{q}S$. It follows that  the splitting dimension of $S$ is $\delta = \dim(D_S)$ and since $D_S\cap \frac{1}{q}S = D_S\cap \frac{1}{q}M_{D_S}$, we have 
\begin{equation*}
    r_F(R) = \lim_{e \to \infty} \frac{a_e}{p^{e\delta}} = \Vol_{M_{D_S}}(P_{\overline{S}} \cap D_S)
\end{equation*}
as desired.

We now show that  the $F$-splitting prime is the monomial ideal generated by the monomials contained in the complement of $D_S$. 

 Supposing $w \in S$, we need to check that $w \in D_S$ if and only there is some $q = p^e > 0$ and $a \in \frac{1}{q}M$ with $\pi_a \in \Hom_R(F^e_*R, R)$ satisfying $\pi_a(x^{\frac{1}{q}w}) = 1$ or equivalently $a + \frac{1}{q}w = 0$.

 Assume first that $w \in D_S$. Setting $a = - \frac{1}{q} w$, we claim $\pi_a \in \Hom_R(F^e_*R,R)$ for $q = p^e \gg 0$. Indeed, as $w \in S$, for the extremal rays $\rho \prec \sigma$ with primitive vectors  $v_{\rho}$ we have $\< w, v_{\rho}\> \geq 0$. If $q \gg 0$, we also have $\< \frac{1}{q}w, v_{\rho}\> < 1$, whence $\<a , v_{\rho} \> > -1$ for $a = - \frac{1}{q} w$ and Theorem~\ref{thm:maps}~(1) is satisfied. If $D$ is an RUF, we know that $w \in D$ as $w \in D_S$ and so also $a = - \frac{1}{q} w \in \langle D \rangle$ and $a \in \frac{1}{q}M_D$. By Theorem~\ref{thm:maps}~(2)~(a)~(i), it follows that $\pi_a \in \Hom_R(F^e_*R,R)$ with $\pi_a(x^{\frac{1}{q}w}) = 1$ and so $x^w$ is not in the splitting prime.

 For the other direction, suppose that $w \not\in D_S$, i.e., there is some RUF $D$ with $w \not\in D$. Let $q = p^e > 0$ and $a \in \frac{1}{q}M$ with $\pi_a \in \Hom_R(F^e_*R, R)$. If $a + {\frac{1}{q}w} \not\in M$, then $\pi_a(x^{\frac{1}{q}w}) = 0 \neq 1$, so assume now also that $a + {\frac{1}{q}w} \in M$. If $a \not\in \langle D \rangle$, then by Theorem~\ref{thm:maps}~(2)~(b) there exists $v \in D^*$ with $\langle a, v \rangle > 0$. Since $\langle w, v \rangle \geq 0$ as $w \in S$, we have $\langle a + {\frac{1}{q}w}, v \rangle > 0$ and so $a + {\frac{1}{q}w} \neq 0$ and so $\pi_a(x^{\frac{1}{q}w}) \neq 1$. Else, if $a \in \langle D \rangle$, we have that $a + {\frac{1}{q}w} \not\in D$ as $w \not\in D$ and again $a + {\frac{1}{q}w} \neq 0$ so $\pi_a(x^{\frac{1}{q}w}) \neq 1$.
\end{proof}

\begin{proposition}\label{lem:splittingratioSignature}
For a seminormal monoid ring $R$ we have $r_F(R)\leq s(\overline{R})$. 
\end{proposition}

\begin{proof}
Note that 
\begin{multline*} r_F(R) = \Vol_{M_{D_S}}(P_S) = \Vol_{M_{D_S}}(P_{\overline{S}}\cap D_S)  = \frac{\Vol_{M}(P_{\overline{S}}\cap D_S)}{[M\cap \<D_S\>:M_{D_S}]} \\
\leq \frac{\Vol_{M}(P_{\overline{S}})}{[M\cap \<D_S\>:M_{D_S}]} = \frac{s(\overline{R})}{[M\cap \<D_S\>:M_{D_S}]}.\end{multline*}

The inequality follows from the fact that the lattice volume of a face of a polytope is bounded by the lattice volume of the polytope. Indeed, if $F$ is a facet of a lattice polytope $P$, then taking a lattice point $p\in P\smallsetminus F$ and letting $Q = \mathrm{conv} (F,p)\subseteq P$, the lattice volume of $Q$ is given by $\Vol_M(Q) = \Vol_M(F) \cdot h$, where $h$ is the distance from $p$ to $F$. Since $\Vol_M(P)\geq \Vol_M(Q)$, the claim follows for facets. For arbitrary faces we proceed by induction. The last equality follows from \cite[Theorem 5.1]{WatanabeYoshidaMinimalRelative}.
\end{proof}

The following example shows that for a seminormal not normal ring the $F$-splitting ratio can be $1$, in particular, equality can be obtained in Proposition \ref{lem:splittingratioSignature}.

\begin{example}
Let $\{e_1,e_2,e_3\}$ be the standard basis of $\Z^3$ and let $C=\<e_1,e_2,e_3\>$. Let $S$ be the seminormal semigroup whose RUFs are $D_1=\<e_2,e_3\>$ and $D_2 =\<e_1,e_3\>$ with $M_{D_1}=\bZ 2e_2 +\bZ e_3$, and $M_{D_2}= \bZ 2e_1 +\bZ e_3$. Then $D_S=\<e_3\>$ and $M_{D_S} = \mathbf{Z} e_3$. We have that $P_{\overline{S}}=\mathrm{conv}\<0,e_1,e_2,e_3, e_1+e_2, e_1+e_3, e_2+e_3, e_1+e_2+e_3\>$ is the standard cube, so $P_S=\mathrm{conv}\<0,e_3\>$, which has volume $1$ with respect to $M_{D_S}$.
\end{example}

\begin{remark}
Recently Jeffries, Nakajima, Smirnov, and Watanabe have shown that the only values of the $F$-signature of a normal monoid ring that are greater than $1/2$ can be $1,2/3,$ and $11/20$, see \cite[Theorem 5.12]{JNSW}. By Proposition \ref{lem:splittingratioSignature}, the same result follows for the $F$-splitting ratio. 
\end{remark}
 
The following example shows that  the $F$-splitting ratio of a ring $R$ and the $F$-signature of its normalisation can be arbitrarily far away from each other within the bounds  $0<s(\overline{R}),r_F(R)\leq 1$.

\begin{example}
Let $S_n$ be the seminormal monoid whose normalization is the monoid corresponding to the cone $C=\<(1,0),(0,1)\>$ 
(i.e.,  $k[x,y]$), where $M_{\<(1,0)\>} = n\Z$ and $M_{\<(0,1)\>} = \Z$. Then 
$P_{\overline{S}}$ is the unit square, so $s(k[\overline{S}])=1$. On the other hand, $D_S = \<(1,0)\>$ and thus $r(k[S_n]) = \Vol_{n\Z} ([0,1])= \frac{1}{n}$.
\end{example}

\begin{remark}
    With Theorem \ref{thm:Fsplittinratio}, it is straightforward to give examples of $F$-split monoid algebras with $F$-splitting primes of arbitrary dimensions. For example, let $S$ be the seminormal monoid whose normalization is the monoid corresponding to the cone generated by the standard basis $e_1, \ldots, e_n$ of $\mathbb{Z}^n$ but where $M_{\langle e_1, \ldots, e_d\rangle} = me_1 + \cdots + me_d$ for some positive integer $m$ relatively prime to $p$. Then $R = k[S]$ has dimension $n$, and the $F$-splitting prime has dimension $d$.
\end{remark}

\section{Test ideals of seminormal monoid algebras}\label{sec:testideals}

The goal of this section is to compute the test ideals of seminormal monoid algebras.

\begin{lemma}
 When $R = k[S]$ is an  monoid algebra over a perfect field $k$, every $\sC$-fixed ideal is monomial.
\end{lemma}

\begin{proof}
 Fix $0 \neq g \in R$. Let us first check that, for all sufficiently large $q = p^e$ (depending on $g$), $\sum_{\phi \in \Hom_R(F^e_*R,R)} \phi(F^e_* gR)$ is a monomial ideal.  Write
 \[
 g = \alpha_1 x^{u_1} + \cdots +\alpha_n x^{u_n}
 \]
 where $\alpha_1, \ldots, \alpha_n \in k^{\times}$ and  $u_1, \ldots, u_n \in S$, with $u_i \neq u_j$ for $i \neq j$.  Take $q = p^e \gg 0$ sufficiently large that $u_i - u_j \not\in q M$ for all $i \neq j$.  Since $\Hom_R(F^e_*R,R)$ is generated by $\pi_\alpha$ for $\alpha \in \frac{1}{q}M$ satisfying the conditions of \autoref{thm:maps}, and $F^e_*R$ is generated as an $R$-module by $x^u$ for $u \in \frac{1}{q}S$, it suffices to check that 
 \[
 \pi_a(g^{1/q}x^u) = \alpha_1^{1/q} \pi_a(x^{u+\frac{1}{q}u_1}) + \cdots + \alpha_n^{1/q} \pi_a(x^{u+\frac{1}{q}u_n})
 \]
 is a monomial for each $a$ and $u$.  But each $\pi_a(x^{u + \frac{1}{q}u_i})$ is monomial, and is nonzero if and only if $a + u + \frac{1}{q}u_i \in M$.  If both $\pi_a(x^{u + \frac{1}{q}u_i})$ and $\pi_a(x^{u + \frac{1}{q}u_j})$ are nonzero for $i \neq j$, then it follows that $\frac{1}{q}u_i - \frac{1}{q}u_j \in M$, contradicting that $q$ was chosen so that $u_i - u_j \not\in qM$.  Thus, we have that at most one of the terms $\pi_a(x^{u + \frac{1}{q}u_i})$ is nonzero, so that  $\pi_a(g^{1/q}x^u)$ is monomial.  
 
 Suppose now that $I = \langle g_1, \ldots, g_N \rangle$ is a $\sC$-fixed ideal.  Iterating \eqref{eq.fixed}, it is easily shown that 
 \[
 I = \sum_{e>e_0} \sum_{\phi\in \Hom_R(F_*^eR,R)} \phi (F^e_* I) = \sum_{i=1}^N \sum_{e>e_0} \sum_{\phi\in \Hom_R(F_*^eR,R)} \phi (F^e_* g_i R)
 \]
 for any $e_0 > 0$.  

 We choose $e_0 \gg 0$ so that the argument above works for $g_1, \ldots, g_N$. A standard argument shows that one can replace $e>e_0$ by $e>0$ in the sum above. It follows that $I$ is a sum of monomial ideals.
\end{proof}

If $R = k[S]$ is a seminormal monoid algebra over a perfect field $k$ and $D \prec C$ is a face, we put $I_D = \langle x^u \mid u \in S \setminus D \rangle$.  If $I$ is any monomial ideal of $R$, then we say that $I$ intersects a face $D$ if and only if $I \not\subseteq I_D$. Note that $I_D$ is a monomial prime ideal of $R$, and if $D, D' \prec C$ then $I_{D} + I_{D'} = I_{D \cap D'}$.  If $D_1, \ldots, D_n \prec C$ are faces, put 
\[I_{D_1 \cup D_2 \cup \cdots \cup D_n} = \langle x^u \mid u \in S \setminus (D_1 \cup D_2 \cup \cdots \cup D_n)  \rangle
= I_{D_1} \cap \cdots \cap I_{D_n},
\]
which is a radical monomial ideal of $R$.  Our goal is to prove the following result, which  characterizes the $\sC$-fixed ideals of $R$.

\begin{theorem}
\label{thm.fixedideals}
 If $R = k[S]$ is a seminormal  monoid algebra over a perfect field $k$, the $\sC$-fixed ideals of $R$ are precisely those of the form
 \[
 I_{D_1 \cup D_2 \cup \cdots \cup D_n}
 \]
for faces $D_1, \ldots, D_n \prec C$ satisfying the following conditions.
 \begin{enumerate}
 \item
 Each $D_i$ is an intersection of proper RUFs of $C$.
 \item
 Every maximal $\mbox{pRUF}$ appears as one of the  $D_i$.
 \end{enumerate}
  \end{theorem}

\noindent
We have largely broken the proof into the series of lemmas below.

\begin{lemma}\label{lem:interior}
 Suppose $I$ is a $\sC$-fixed monomial ideal and $D \prec C$ is a face of $C$ not contained in any $p$-face. If $I$ intersects $D$, then $x^w \in I$ for all $w \in S \cap \interior(D)$.
\end{lemma}

\begin{proof}
 Let $x^u \in I$ with $u \in S \cap D$.  Replacing with an appropriate multiple in $k[S\cap \interior(D)]$ if necessary, we may assume that $u \in \interior(D)$. Let $w\in S\cap \interior(D)$.
 For $q = p^e \gg 0$ sufficiently large, we have as well that $w - \frac{2}{q} u \in \interior(D) \cap \frac{1}{q}M_D$.  Put $b = qw - 2u \in M_D \cap \interior(D)$, so that $b \in S$ and $x^{b + u} \in I$.

Letting $a = \frac{1}{q}u$, we check that $a$ satisfies the conditions of  \autoref{thm:maps}. (1) and (2) (b) follow immediately from the fact that $a\in C$. 
For (2) (a), let $\widetilde{D}$ be a RUF not properly contained in a $p$-face and assume $a\in \langle \widetilde{D}\rangle$. Then $\widetilde{D}$ is not a $p$-face, since $a\in \interior{D}$ and this would imply $D$ is contained in a $p$-face, a contradiction. So (2) (ii) does not occur. If $\widetilde{D}$ is not a $p$-face, then  since $u\in S\cap \widetilde{D}$, we have $u\in M_{\widetilde{D}}$, so $a\in \frac{1}{q}M_{\widetilde{D}}$ and (2) (i) follows. 

Thus, $\pi_a \in \Hom_R(F^e_*R,R)$, whence $\pi_a(F^e_*x^{b+u}) = x^w \in I$ as we have

\[
a+\frac{1}{q}(b+u) = \frac{1}{q}(b+2u)=\frac{1}{q}(qw) = w \in M,
\]
$x^{b + u} \in I$, and $I$ is $\sC$-fixed.
\end{proof}

\begin{lemma}
\label{lem.canthit}
 If $D \prec C$ is a $p$-face and $u \in S \cap \langle D \rangle$, then 
 $$x^u \not\in \sum_{e > 0} \sum_{\phi \in \Hom_R(F^e_*R,R)}\phi(F^e_*R).$$
 In particular, $\sC$-fixed ideals cannot intersect $p$-faces. 
\end{lemma}

\begin{proof}
 By way of contradiction, suppose there is $q=p^e$, $a \in \frac{1}{q}M$ satisfying  the conditions of \autoref{thm:maps}, and  $w \in \frac{1}{q}S$ such that  $\pi_a(x^{w}) = x^u$.  Thus, we have $a + w = u \in M_D$.  This implies $a,w \in D$. By \ref{thm:maps} 2) a) ii), we have $a \not\in \frac{1}{q}\widetilde{M}_{D}$ where $\widetilde M_D = \{ x \in M \cap \< D \> \, | \, mx \in M_D \mbox{ for some } m \in \Z \mbox{ with } \gcd(m,p) = 1 \}$. On the other hand, $a=u-w\in \frac{1}{q}S\subseteq \frac{1}{q}M_D\subseteq \frac{1}{q}M_{\widetilde{D}}$, a contradiction. 
\end{proof}

\begin{lemma}
\label{lem.unionRUFs}
 Let $I$ be a $\sC$-fixed monomial ideal. Let $u \in S$ and  
 $$D = \bigcap_{ \substack{\mbox{ $ u \in \widetilde{D}$ } \\ \mbox{$\widetilde D$ a RUF}} } \widetilde D,$$ 
 the intersection of all proper RUFs containing $u$. If $D$ intersects $I$, then $x^u \in I$.  
\end{lemma}

\begin{proof}
By Lemma \ref{lem.canthit} $I$ does not intersect any $p$-face, so $D$ is not contained in any $p$-face and the assumption of Lemma \ref{lem:interior} are satisfied.  Thus we have that $\interior(D) \cap S \subseteq I$.  Let $D'$ be the smallest face containing $u$, so that $u \in \interior(D') \cap S$. As $u \in D$, we must have $D' \prec D$, and so $D^* \prec D'^*$.  Pick $-b \in \interior(D) \cap M_D$.   Choose an integer $\ell \gg 0$ sufficiently large so that $\ell \langle u, v_\rho \rangle > \langle -b, v_\rho \rangle$ for all primitive vectors $v_\rho \not\in D'^*$.  Setting $w = b + \ell u$, for any primitive vector $v_\rho$ we then have:
\begin{itemize}
 \item
 $\langle w, v_\rho \rangle > 0$ if $v_\rho \not\in D'^*$;
 \item
 $\langle w, v_\rho \rangle < 0$ if $v_\rho \in D'^* \setminus D^*$;
 \item
 and $w \in \langle D \rangle$, so that $\langle w, v_\rho \rangle = 0$ for all $v_\rho \in D^*$.
\end{itemize}
Choose $q = p^e \gg 0$ so that $q > \ell$ and $-q < \langle w,v_\rho \rangle$ for all $v_\rho \in D'^* \setminus D^*$, and put $a = \frac{1}{q}w \in \frac{1}{q}M$.  

We now check that $a$ satisfies the conditions of \autoref{thm:maps}.  Condition (1) follows from our choice of $q$, so consider  $\widetilde D \prec C$ a RUF.    Suppose first  $a \not \in \langle \widetilde D \rangle$, then so too $w = qa \not \in  \langle \widetilde D \rangle$.  If $D' \prec \widetilde D$, then we would have $u \in D' \prec \widetilde D$ and so $D \prec \widetilde D$, which would contradict $w \not \in \langle \widetilde D \rangle$ as $w \in \langle D \rangle$.  Thus, we must have $D' \not \prec \widetilde D$, so $\widetilde D^* \not \prec D'^*$.  Taking $v_\rho \in \widetilde D^* \setminus D'^*$, we have $\langle a, v_\rho \rangle > 0$ so that condition (2a) or (3a) is satisfied, respectively.  Else, suppose instead that $a \in \langle \widetilde D \rangle$ and so $w \in \langle \widetilde D \rangle$.  Thus, for all $v_\rho \in \widetilde D^*$, we have $\langle w, v_\rho \rangle = 0$, and so necessarily $v_\rho \in D^*$ as well.  It follows that $\widetilde D^* \prec D^*$, so $D \prec \widetilde D$.  Since $D$ intersects $I$, we know that $D$ is not contained in any $p$-face.  Thus, $\widetilde D$ cannot be a maximal pRUF, and we may assume $\widetilde D$ is not contained in any $p$-face.  We have $w \in M_D \subseteq M_{\widetilde D}$, so that $a \in \frac{1}{q}M_{\widetilde D}$ and (2a) is satisfied.  Thus, it follows $\pi_a \in \Hom_R(F^e_*R,R)$.

Let $y = -b + (q - \ell) u \in \interior(D) \cap M_D \subseteq S$ so that $x^y \in I$. Then $\pi_a(x^{\frac{1}{q}y}) = \frac{1}{q}(w + y) = \frac{1}{q}(b + \ell u -b + (q -\ell) u) = u$ as desired.
\end{proof}

\begin{lemma}
\label{lem.intheface}
Suppose $D \prec C$ is a RUF and $0 \neq \pi_a(x^u) \in D$ for some $a \in \frac{1}{q}M$, $\pi_a \in \Hom_R(F^e_*R,R)$, and $u \in \frac{1}{q}S$.  Then we must have $a \in \langle D \rangle$ and $u \in  D$.
\end{lemma}

\begin{proof}
 Let $\pi_a(x^u) = x^w \in R$, so that $a + u = w \in S \cap D$.  If $a \not\in \langle D \rangle$, then we know $\langle a, v_\rho \rangle > 0$ for some primitive vector $v_\rho \in D^*$.  Since $u \in \frac{1}{q}S$, also $\langle u, v_\rho \rangle \geq 0$ and it follows $\langle w, v_\rho \rangle > 0$ contradicting $w \in D$.  Thus, it follows we must have $a \in \langle D \rangle$ and so also $u = w - u \in \langle D \rangle \cap \frac{1}{q}S \subseteq D$.
\end{proof}

\begin{proof}[Proof of \autoref{thm.fixedideals}.]  Consider first an ideal
\[
 I = I_{D_1 \cup D_2 \cup \cdots \cup D_n}
 \]
for faces $D_1, \ldots, D_n \prec C$ where
each $D_i$ is an intersection of  RUFs $C$, and also 
 every maximal $\mbox{pRUF}$ appears as one of the  $D_i$.  We need to show that $I$ is $\sC$-fixed.
 
 If $a \in \frac{1}{q}M$ and $\pi_a \in \Hom_R(F^e_*R,R)$, we first argue that $\pi_a(F^e_*I_{D_i}) \subseteq I_{D_i}$ for each $i$.  If $\pi_a(x^u) \in D_i$ for some $u \in \frac{1}{q}S$, we know by \ref{lem.intheface} that $u \in \widetilde D$ for any RUF $\widetilde D$ containing $D$.  Since $D_i$ is an intersection of RUFs, it follows that $\pi_a(x^u) \in D_i$ for some $u \in \frac{1}{q}S$ implies $u \in D_i$.  Conversely, $u \not \in D_i$ gives that $\pi_a(x^u) \not \in D_i$ and it follows $\pi_a(F^e_*I_{D_i}) \subseteq I_{D_i}$.  This condition is compatible with intersections, so we conclude that $\pi_a(F^e_*I) \subseteq I$ and thus $\phi(F^e_*I) \subseteq I$ for any $\phi \in \Hom_R(F^e_* R, R)$ and $e \in \N$. It follows that $$\sum_{e > 0} \sum_{\phi \in \Hom_R(F^e_*R,R)}\phi(F^e_*R) \subseteq I.$$ 
 
Suppose now that $x^u \in I$.  Since the maximal pRUF's are among the $D_i$'s, we know that $u$ is not in any $p$-face.  If $a = \frac{1}{p} u$, then it is easy to see from \autoref{thm:maps} that $\pi_a \in \Hom_R(F^e_*R,R)$.  We have that $\pi_a(x^{\frac{(p-1)u}{p}}) = u$ and $x^{(p-1)u}\in I$.  Thus, it follows that $I$ is $\sC$-fixed.

Conversely, if $I$ is a $\sC$-fixed ideal, it follows from \autoref{lem.unionRUFs} that $\{ u \in S \mid x^u \not \in I \}$ is a union of lattice points in intersections of RUFs.  Moreover, by \autoref{lem.canthit}, we also know this set contains every lattice point in every maximal pRUF.  Thus, $I$ must be of the desired form.
\end{proof}

\begin{proof}[Proof of Corollary \ref{cor:Cfixed}] 
     This follows immediately from Theorem \ref{thm.fixedideals}. The first part follows from the definition of $I_{D_1 \cup \cdots \cup D_n}$.
    For the second part, to get the the test ideal which  is the smallest non-zero $\sC$-fixed ideal, we take  $I_{D_1 \cup \cdots \cup D_n}$, where $D_1, \ldots, D_n$ runs over all RUFs. 
  For the third part, to get the non-$F$-pure ideal which is the largest $\sC$ fixed ideal, we take the minimal number of faces away, and get 
  $\sigma(R) = I_{D_1 \cup \cdots \cup D_n}$ where the $D_i$ run through the $p$-RUFs. 
  \end{proof}

\bibliography{mybib,MainBib}
\end{document}